\documentclass[11pt]{article}
\usepackage{amssymb,latexsym}
\usepackage{epsfig}
\usepackage{eufrak}
\usepackage{amsmath}
\usepackage{mathrsfs}
\usepackage{color}

\newtheorem{theorem}{Theorem}
\newtheorem{lemma}{Lemma}

\newcommand{\be}{\begin{equation}}
\newcommand{\ee}{\end{equation}}
\newcommand{\bee}{\begin{eqnarray*}}
\newcommand{\eee}{\end{eqnarray*}}
\newcommand{\bel}{\begin{eqnarray}}
\newcommand{\eel}{\end{eqnarray}}
\newcommand{\bec}{\begin{cases}}
\newcommand{\eec}{\end{cases}}
\newcommand{\bem}{\begin{bmatrix}}
\newcommand{\eem}{\end{bmatrix}}

\newcommand{\bed}{\begin{description}}
\newcommand{\eed}{\end{description}}
\newcommand{\bei}{\begin{itemize}}
\newcommand{\eei}{\end{itemize}}
\newcommand{\ben}{\begin{enumerate}}
\newcommand{\een}{\end{enumerate}}

\newcommand{\beL}{\begin{lemma}}
\newcommand{\eeL}{\end{lemma}}
\newcommand{\beT}{\begin{theorem}}
\newcommand{\eeT}{\end{theorem}}
\newcommand{\sect}{\section}

\newcommand{\bpf}{\begin{pf}}
\newcommand{\epf}{\end{pf}}

\newcommand{\pfbox}{\hfill\mbox{$\Box$}}

\newenvironment{pf}{\paragraph*{Proof{\rm.}}}{\pfbox\bigskip}

\begin{document}

\title{\bf Fast Universal Algorithms for Robustness Analysis\thanks{This
research was supported in part by grants from NASA (NCC5-573) and
LEQSF (NASA /LEQSF(2001-04)-01).} }

\author{Xinjia Chen, Kemin Zhou, and Jorge L. Aravena\\
Department of Electrical and Computer Engineering\\
Louisiana State University\\
Baton Rouge, LA 70803\\
\{chan,kemin,aravena\}@ece.lsu.edu\\
Tel: (225)578-\{5533,5537\}\\
Fax: (225) 578-5200 }

\date{March 2003}

\maketitle

\begin{abstract}

In this paper, we develop efficient randomized algorithms for estimating
probabilistic robustness margin and constructing robustness degradation curve for
uncertain dynamic systems.  One remarkable feature of these algorithms is their
universal applicability to robustness analysis problems with
arbitrary robustness requirements and uncertainty bounding set.
In contrast to existing probabilistic methods, our approach does not
depend on the feasibility of computing deterministic robustness margin.  We have developed efficient methods
such as probabilistic comparison, probabilistic bisection and
backward iteration to facilitate the computation.  In particular, confidence
interval
for binomial random variables has been frequently used
in the estimation of probabilistic robustness margin
and in the accuracy evaluation of estimating robustness degradation
function.
Motivated by the importance of fast computing of binomial confidence
interval
in the context of probabilistic robustness analysis, we have derived an
explicit formula
for constructing the confidence
interval of binomial parameter with guaranteed coverage probability.
The formula overcomes the limitation of
normal approximation which is asymptotic
in nature and thus inevitably introduce
unknown errors in applications.
Moreover, the formula is extremely simple and very tight in
comparison with classic Clopper-Pearson's
approach.

\end{abstract}

\section{Introduction}

In recent years, there have been growing interest
on the development of probabilistic methods for robustness analysis and
design problems aimed at overcoming the computational complexity
and the issue of conservatism of deterministic
worst case framework \cite{RS, SR, MS, KT, TD, BLT, BL, BP, Polyak, CDT1,
CDT2, C1, C2, C3, vind, vind2, SB}.  In the deterministic worst case framework,
one is interested in computing the {\it deterministic robustness margin}
which is defined as the maximal radius of uncertainty set such that the robustness requirement is
guaranteed everywhere.  However, it should be borne in mind that the uncertainty set
may include worst cases which never happen in reality.
Instead of seeking the worst case guarantee, it is sometimes ``acceptable''
that the robustness requirement is satisfied for most of the cases.
It has been demonstrated that the proportion of systems
guaranteeing  the robustness requirement
can be close to $1$ even if the radii of uncertainty set are
much larger than the deterministic robustness margin.
The idea of sacrificing extreme cases of uncertainty has become a new paradigm of robust control.
In particular, the concept of {\it probabilistic robustness margin}
has been recently established \cite{BLT,BL,TD,BLT2, BP, CDT1, Polyak}.  Moreover,
to provide more insight for the robustness of the uncertain system, the concept
of {\it robustness degradation function} has been developed by a number of researchers \cite{BLT, CDT1}.
For example,  Barmish and Lagoa \cite{BL} have
constructed a curve of robustness margin amplification versus risk in a probabilistic setting.
In a similar spirit, Calafiore, Dabbene and Tempo \cite{CDT1, CDT2} have
constructed a probability degradation function in the context of
real and complex parametric uncertainty.
Such a function describes quantitatively the relationship between
the proportion of systems
guaranteeing  the robustness requirement and the radius of uncertainty set.
Clearly, it can serve as a guide for control engineers in
evaluating the robustness of
a control system once a controller design is completed.  It is important to
note that the robustness degradation function and the probabilistic robustness margin can be estimated in a distribution
free manner.  This can be justified by the Truncation Theory established by Barmish, Lagoa
and Tempo \cite{BLT} and can also be illustrated by relaxing the deterministic worst-case paradigm.

Existing techniques for constructing the robustness degradation
curve relies on the feasibility of computing the deterministic robustness margin.
Obtaining the probabilistic robustness margin is possible only when the
robustness degradation curve is constructed.
However, the computation of deterministic robustness
margin is possible only when the robustness requirement ${\bf P}$ is simple
and the uncertainty bounding set is of some special structure.
For example, ${\bf P}$ is stability requirement and uncertainty bounding sets
are spectral normed balls.
In that case, deterministic analysis methods such as $\mu$-theory
can be applied to compute the deterministic robustness margin.  The
deterministic robustness margin is then taken
as a starting point of uncertainty radius interval for which the robustness
degradation curve is to be constructed.
In general, the problem of computing the deterministic robustness margin is
not tractable.
Hence, to construct a robustness degradation curve of practical interest,
the selection of
uncertainty radius interval is itself a question.  Clearly, the range of
uncertainty radius for which robustness degradation
curve is significantly below $1$ is not of practical interest since only a
small risk can be tolerated in reality.
From application point of view, it is only needed to construct robustness
degradation curve for the range of
uncertainty radius such that the curve is above an
{\it a priori} specified level $1-\epsilon$ where risk parameter $\epsilon
\in (0,1)$ is acceptably small.

In this work, we consider robustness analysis
problems with arbitrary robustness requirement and uncertainty bounding set.
We develop efficient randomized algorithms for estimating probabilistic
robustness margin $\rho_{\epsilon}$
which is defined as the maximal uncertainty radius
such that the probability of guaranteeing the robust requirements
is at least $1-\epsilon$.  We have also developed fast algorithms for
constructing robustness degradation curve which is above a priori specified
level $1-\epsilon$.
In particular, we have developed efficient mechanisms such as probabilistic
comparison,
probabilistic bisection and backward iteration to reduce the computational
complexity.  Complexity of probabilistic comparison
techniques are rigorously quantified.

In our algorithms, confidence interval
for binomial random variables has been frequently used
to improve the efficiency of estimating probabilistic robustness margin
and in the accuracy evaluation of robustness degradation function.
Obviously, fast construction of binomial confidence interval is
important to the efficiency of the randomized algorithm.
Therefore, we have derived an explicit formula
for constructing the confidence
interval of binomial parameter with guaranteed coverage probability \cite{C4}.
The formula overcomes the limitation of
normal approximation which is asymptotic
in nature and thus inevitably introduce
unknown errors in applications.
Moreover, the formula is extremely simple and very tight in
comparison with classic Clopper-Pearson's
approach.

The paper is organized as follows.  Section 2 is the problem formulation.
Section 2 discusses binomial confidence interval. Section 3 is devoted to
probabilistic robustness margin.
Section 4 presents algorithms for constructing robustness degradation curve.
Illustrative examples are given in Section 5.  Section 6 is the conclusion.  Proofs are included in Appendix.

\section{Problem Formulations}

We adopt the assumption, from the classical robust control framework, that the
uncertainty is deterministic and bounded.
We formulate a general robustness analysis problem  in a similar way as \cite{C3} as follows.

Let ${\bf P}$ denote a robustness requirement.
The definition of ${\bf P}$ can be a fairly complicated
combination of the following:

\begin{itemize}

\item Stability or ${\cal D}$-stability;

\item $H_{\infty}$ norm of the closed loop transfer function;

\item  Time specifications such as
overshoot, rise time, settling time and steady state error.

\end{itemize}

Let ${\cal B}(r)$ denote the
set of uncertainties with size smaller than $r$. In applications, we are usually dealing with
uncertainty sets such as the following:
\begin{itemize}
\item $l_p$ ball
$
{\cal B}^p(r) :=\left\{ \Delta \in {\bf R}^n:\;||\Delta||_p \leq r \right\}
$
where $||.|_p$ denotes the $l_p$ norm and $p = 1, 2, \cdots,\infty$.
In particular, ${\cal B}^{\infty} (r)$ denotes a box.
\item  Spectral norm ball
$
{\cal B}_{\sigma}(r) := \{ \Delta \in {\bf \Delta}: \bar{\sigma}(\Delta) \leq r \}$
where $\bar{\sigma}(\Delta)$ denotes the largest singular value of $\Delta$.  The
 class of allowable perturbations is
\be
\label{str}
{\bf \Delta} : =\{ {\rm blockdiag} [q_1I_{r_1}, \cdots, q_sI_{r_s}, \Delta_1, \cdots, \Delta_c  ]\}
\ee
where $q_i \in \Bbb{F}, \;i=1, \cdots,s$ are scalar parameters with multiplicity $r_1,\cdots,r_s$ and
$\Delta_i \in \Bbb{F}^{n_i \times m_i}, \;i = 1, \cdots, c$ are possibly repeated full blocks.
Here $\Bbb{F}$ is either the complex field ${\bf C}$ or the real field ${\bf R}$.

\item Homogeneous star-shaped bounding set
$
{\cal B}_H (r):=\left\{ r (\Delta - \Delta_0) + \Delta_0 :
 \Delta \in Q \right\}
 $
 where $Q \subset {\bf R}^n$ and $\Delta_0 \in Q$
 (see \cite{BLT} for a detailed illustration).
\end{itemize}
Throughout this paper, ${\cal B} (r)$ refers to any type of uncertainty set described above.
Define a function $\ell ( . )$  such that,
 for any $X$,
 \[
 \ell (X) := \min \{ r : X \in {\cal B} (r) \},
 \]
  i.e.,
${\cal B}(\ell (X))$ includes $X$ exactly in the boundary.  By such definition,
\[
\ell (X) = \min \left\{ r :  \frac{X- \Delta_0}{r} + \Delta_0 \in Q \right\},
\]
\[
\ell (X) = \bar{\sigma}(X),
\]
and
\[
\ell (X) = ||X||_p
\]
in the context of homogeneous star-shaped bounding set, spectral norm ball and $l_p$ ball respectively.

To allow the robustness analysis be performed in a distribution-free manner,
we introduce the notion of {\it proportion} as follows.  For any $\Delta \in
{\cal B}(r)$
there is an associated system $G(\Delta)$.  Define {\it proportion}
\[
\Bbb{P} (r) := \frac{{\rm vol} ( \{\Delta \in {\cal B}(r):
{\rm The \; associated \; system}
\; G(\Delta) \; {\rm guarantees} \; {\bf P} \}   ) } { {\rm vol} ( {\cal
B}(r) ) }
\]
with
\[
{\rm vol} (S) := \int_{q \in S} dq,
\]
where the notion of $dq$ is illustrated as follows:

\begin{itemize}

\item (I): If $q =[x_{rs}]_{n \times m}$ is a real matrix in ${\bf R}^{n
\times m}$,
then $dq = \prod_{r=1}^n \prod_{s=1}^m dx_{rs}$.

\item (II): If $q =[x_{rs} + j y_{rs}]_{n \times m}$ is a complex
matrix in ${\bf C}^{n \times m}$, then $dq = \prod_{r=1}^n \prod_{s=1}^m (
dx_{rs} dy_{rs} )$.

\item (III): If $q \in {\bf \Delta}$, i.e., $q$ possess a block structure
defined by (\ref{str}), then
$dq = (\prod_{i=1}^s d q_i) (\prod _{i=1}^c d \Delta_i)$ where the notion
of $dq_i$ and $d \Delta_i$ is defined by (I) and (II).

\end{itemize}

It follows that $\Bbb{P} (r)$ is a reasonable measure of the robustness of
the system \cite{CDT1,TD2}.
In the worst case deterministic framework,
we are only interested in knowing if ${\bf P}$ is guaranteed for every
$\Delta$.
However, one should bear in mind that the uncertainty set in
our model may include worst cases which never happen in reality.
Thus, it would be ``acceptable'' in many applications
if the robustness requirement ${\bf P}$ is satisfied for most of the cases.
Hence,
due to the inaccuracy of the model, we should also
obtain the value of $\Bbb{P} (r)$ for uncertainty radius $r$ which exceeds
the deterministic robustness margin.

Clearly, $\Bbb{P} (r)$ is deterministic
in nature.  However, we can resort to
a probabilistic approach to evaluate $\Bbb{P} (r)$.  To see this, one needs
to
observe that a random variable with
{\it uniform distribution} over ${\cal B} (r)$, denoted by ${\bf \Delta}^u$,
guarantees that
\[
\Pr \{ {\bf \Delta}^u \in S \} = \frac{{\rm vol} (S \bigcap {\cal B}(r)) } {
{\rm vol} ({\cal B}(r))  }
\]
for any $S$, and thus
\[
\Bbb{P} (r) =  \Pr \{
{\rm The \; associated \; system}
\; G({\bf \Delta}^u) \; {\rm guarantees} \; {\bf P} \}.
\]
Define a Bernoulli random variable $X$ such that $X$ takes value $1$
if the associated system $G({\bf \Delta}^u)$ guarantees ${\bf P}$ and takes
value $0$ otherwise.  Then
estimating $\Bbb{P}(r)$ is equivalent to estimating binomial parameter
$\Bbb{P}_X : = {\rm Pr} \{X=1\} = \Bbb{P}(r)$.
It follows that a Monte Carlo method can be employed to estimate $\Bbb{P}
(r)$
based on i.i.d. observations of $X$.

Obviously, the robustness problem will be completely solved if we can
efficiently estimate $\Bbb{P}(r)$ for all $r \in (0, \infty)$.
However, this is infeasible from computational perspective.  In practice,
only a small risk $\epsilon$
can be tolerated by a
system.  Therefore, what is really
important to know is the value of  $\Bbb{P}(r)$ {\it over the range of
uncertain radius $r$ for which
$\Bbb{P}(r)$ is at least $1 - \epsilon$}
where $\epsilon \in (0,1)$ is referred as the {\it risk} parameter in this
paper.  When the {\it deterministic robustness margin}
\[
\rho_0 := \sup \left\{ r: \; G(\Delta) \; {\rm guarantees}
\; {\bf P} \; {\rm for \; any} \; \Delta \in {\cal B}(r) \right\}
\]
can be computed, we can construct the robustness degradation curve as follows:
(I) Estimate proportion for a sequence of discrete values of uncertainty radius which
is started from $r = \rho_0$ and gradually increased;
(II) The construction process is continued until the proportion is below $1 - \epsilon$.
This is the conventional approach and it is feasible only when computing $\rho_0$ is tractable.
Unfortunately, in general there exists no effective
technique for computing the deterministic robustness margin $\rho_0$.
In light of this situation, we establish a new approach
which does not depend on the feasibility of computing the deterministic robustness margin.
Our strategy is to firstly estimate the {\it probabilistic
robustness margin}
\[
\rho(\epsilon):=\sup \{r:\;\Bbb{P}(r) \geq 1 - \epsilon\}
\]
and consequently construct the robust degradation curve in
a backward direction (in which $r$ is decreased) by choosing
the estimate of $\rho(\epsilon)$ as the starting uncertainty radius.
\bigskip

To reduce computational burden, the estimation of probabilistic robustness
margin
relies on the frequent use of binomial confidence interval.  The confidence
interval
is also served as a validation method for the accuracy of estimating
robustness degradation function.
Hence, it is desirable to
quickly construct binomial confidence interval with guaranteed coverage
probability.

\section{Binomial Confidence Intervals}
Clopper and Pearson \cite{Clo} provided a rigorous approach for
constructing binomial confidence interval.
However, the computational complexity involved with this approach is very
high.
The standard technique is to use normal approximation which is not accurate
for rare events.
The coverage probability of the confidence interval derived from normal
approximation can be significantly below the specified confidence level even
for very large sample size.
In the context of robustness analysis, we are dealing with rare events
because the probability
that the robustness requirement is violated is usually very small.  We shall
illustrate these standard methods  as follows.

\subsection{Clopper-Pearson Confidence Limits}
Let the sample size $N$ and confidence parameter $\delta \in (0,1)$ be
fixed.
We refer an observation of $X$
with value 1 as a {\it successful trial}.
Let $K$ denote the number of successful trials during the $N$ i.i.d.
sampling experiments.
Let $k$ be a realization of $K$.  The classic Clopper-Pearson lower
confidence limit $L_{N,k,\delta}$ and
   upper confidence limit $U_{N,k,\delta}$ are given respectively by
\[
L_{N,k,\delta}:=\left\{\begin{array}{ll}
   0 \;\;\;&  {\rm if}\; k=0\\
   \underline{p} \; \;\;\;&
   {\rm if}\; k > 0
\end{array} \right.\;\;\;{\rm and}\;\;\;
U_{N,k,\delta}:=\left\{\begin{array}{ll}
   1  &  {\rm if}\; k=N\\
   \overline{p} \;
   & {\rm if}\; k < N
\end{array} \right.
\]
where $\underline{p} \in (0,1)$ is the solution of the following equation
\be
\label{clp1}
\sum_{j=0}^{k-1} {N \choose j}
   \underline{p}^j (1- \underline{p})^{N-j} = 1-\frac{\delta}{2}
\ee
and $\overline{p} \in
(0,1)$ is the solution of the following equation
\be
\label{clp2}
\sum_{j=0}^{k} {N \choose j}
   \overline{p}^j (1- \overline{p})^{N-j} = \frac{\delta}{2}.
\ee

\subsection{Normal Approximation}

It is easy to see that the equations~(\ref{clp1}) and~(\ref{clp2}) are
hard to solve and thus the confidence limits are difficult
to determine using Clopper-Pearson's approach.
For large sample size, it is computationally
intensive.  To get around the difficulty, normal
approximation has been widely used to develop simple approximate formulas
(see, for example, \cite{Hald} and the references therein).
Let $\Phi(.)$ denote the normal distribution function and
$Z_{\frac{\delta}{2}}$ denote the critical value
such that $\Phi( Z_{\frac{\delta}{2}} ) = 1 - \frac{\delta}{2}$.
It follows from the Central Limit Theorem that, for sufficiently large sample size $N$,
the lower and upper confidence limits can be estimated respectively as
$\widetilde{L} \approx  \frac{k}{N}  - Z_{\frac{\delta}{2}}
\sqrt{ \frac{ \frac{k}{N} (1-\frac{k}{N}) } {N} }$ and $\widetilde{U}
\approx  \frac{k}{N} + Z_{\frac{\delta}{2}} \sqrt{ \frac{ \frac{k}{N}
(1-\frac{k}{N}) } {N} }$.

The critical problem with the normal approximation is that it is of
asymptotic nature.  It is not clear
how large the sample size is sufficient for the approximation error to be
negligible.
Such an asymptotic approach is not good enough for studying the robustness
of control systems.

\subsection{Explicit Formula}

It is desirable to have a simple formula which is rigorous and very
tight for the confidence interval construction. Recently, we have
derived the following simple formula for constructing the confidence
limits \cite{C4} (The proof is provided in Appendix for purpose of
completeness).

\begin{theorem} \label{Tape_Massart}
Let ${\cal L} (N,k,\delta) = \frac{k}{N} + \frac{3}{4} \;
\frac{ 1 - \frac{2k}{N} -
\sqrt{ 1 + 4 \theta \; k ( 1- \frac{k}{N}) } }
{1 + \theta N }$ and ${\cal U} (N,k,\delta) = \frac{k}{N} + \frac{3}{4} \;
\frac{ 1 - \frac{2k}{N} +
\sqrt{ 1 + 4 \theta \; k ( 1- \frac{k}{N}) } } {1 + \theta N }$
with $\theta = \frac{9}{ 8 \ln \frac{2}{\delta} }$.
Then $\Pr \left\{   {\cal L} (N,K,\delta)< \Bbb{P}_X < {\cal U} (N,K,\delta) \right\} > 1 -
\delta$.
\end{theorem}

\begin{figure}[htbp]
\centerline{\psfig{figure=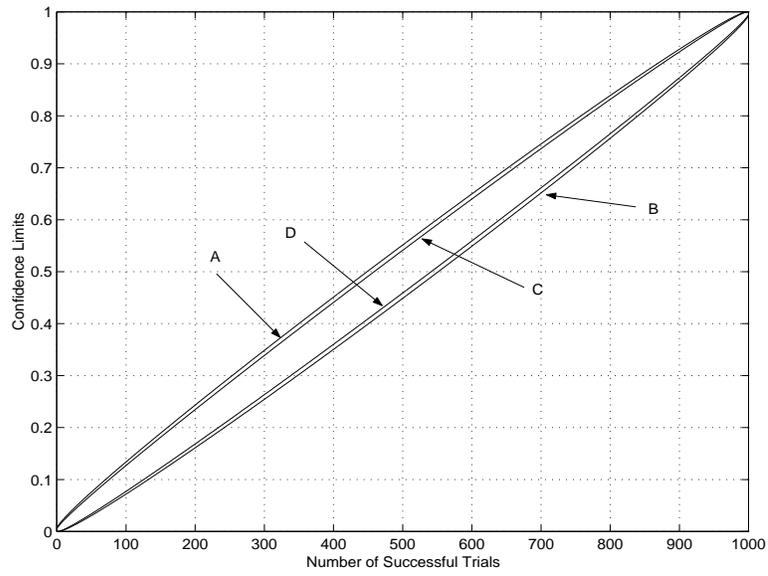, height=3.0in, width=4.0in
}} \caption{ Confidence Interval ($N = 1000, \; \delta = 10^{-2}$) }
\label{fig_1}
\end{figure}

\begin{figure}[htbp]
\centerline{\psfig{figure=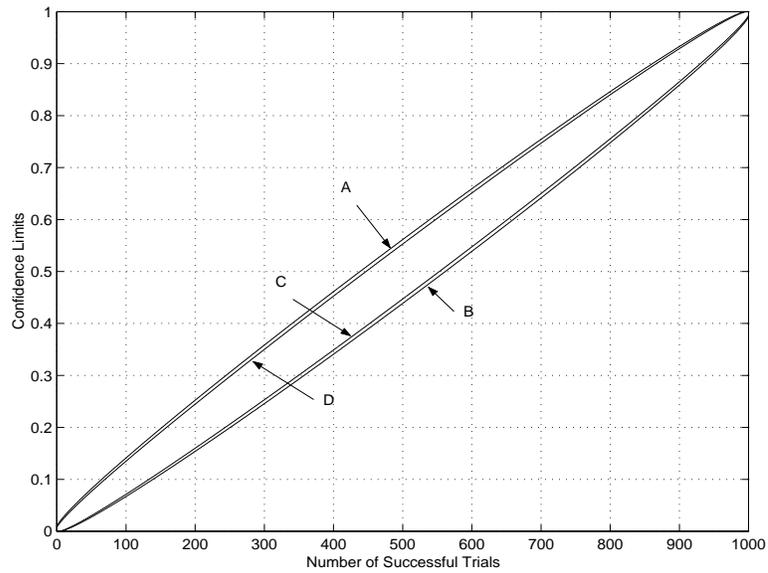, height=3.0in, width=4.0in
}} \caption{ Confidence Interval ($N = 1000, \; \delta = 0.001$) }
\label{fig_2}
\end{figure}

Figures \ref{fig_1} and \ref{fig_2} show the confidence
limits derived by different methods
(curve A and B represent respectively the upper and lower confidence limits computed by Theorem~\ref{Tape_Massart};
curve C and D represent respectively the upper and lower confidence limits calculated by
Clopper-Pearson's method).  It can be seen from these figures
that our formula is very tight in comparison with the Clopper-Pearson's approach.
Obviously, there is no comparison on
the computational complexity.
Our formula is simple enough for hand calculation.
Simplicity of the confidence interval is especially important
in the context of our robustness analysis problem
since the confidence limits are repeatedly used for a large number of
simulations.

\bigskip

\bigskip

\section{Estimating Probabilistic Robustness Margin}

In this section, we shall develop efficient randomized algorithms for
constructing an estimate for $\rho(\epsilon)$.

\subsection{Separable Assumption}

We assume that {\it the robustness degradation curve of the system can be
separated into two parts
by a horizontal line with height $1-\epsilon$}, i.e.,
\[
\Bbb{P}(r) < 1 - \epsilon \; {\rm for \; all} \; r \geq \rho(\epsilon).
\]

\begin{figure}[htbp]
\centerline{\psfig{figure=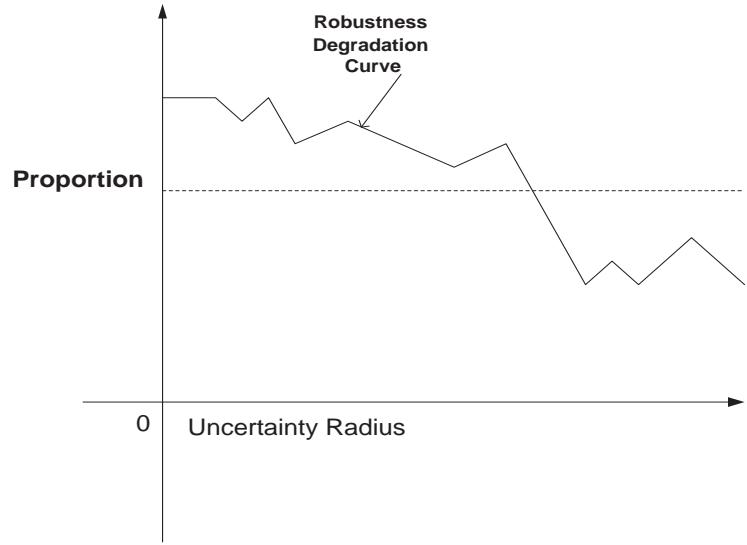, height=3in, width=4.0in
}} \caption{Illustration of Separable Assumption.  The robustness
degradation curve can be
separated as the upper segment and lower segment by the dashed horizontal line
with height $1-\epsilon$.
In this example, the separable assumption is satisfied, while the
non-increasing assumption is violated.}
\label{fig_3}
\end{figure}

We refer such an assumption as the {\it Separable Assumption}.
From an application point of view, this assumption is rather benign.
Our extensive simulation experience indicated that, for small risk parameter
$\epsilon$,
most control systems guarantee the separable assumption.  It should be
noted that it is even much weaker
than assuming that $\Bbb{P}(r)$ is non-decreasing (See illustrative
Figure~\ref{fig_3}).
Moreover, the non-increasing assumption is rather mild.  This can be
explained by a heuristic argument as follows.
Let
\[
{\cal B}^{\bf P} (r) : =  \{\Delta \in {\cal B}(r):
{\rm The \; associated \; system}
\; G(\Delta) \; {\rm guarantees} \; {\bf P} \}.
\]
Then
\[
\Bbb{P}(r) = \frac{{\rm vol} ( {\cal B}^{\bf P} (r) )  } { {\rm vol} ({\cal
B}(r))  }
\]
and
\be
\frac{d \Bbb{P}(r)}{d r} = \frac{ 1 } { {\rm vol} ({\cal B}(r)) }
\left[  \frac{d \; {\rm vol}\{ {\cal B}^{\bf P} (r) \} }{d r} -
\Bbb{P}(r)\;\frac{ d\; {\rm vol} \{ {\cal B} (r) \} }{d r} \right].
\label{der}
\ee
Moreover, due to the constraint of robust
requirement
${\bf P}$, it is true that ${\cal B}^{\bf P} (r + dr) \; \backslash \; {\cal B}^{\bf P} (r)$ is a subset of
${\cal B} (r + dr) \; \backslash \; {\cal B} (r)$ and it follows that
${\rm vol} \{ {\cal B} (r) \}$ increases (as $r$ increases) faster than
${\rm vol}\{ {\cal B}^{\bf P} (r) \}$, i.e., $\frac{d \; {\rm vol}\{ {\cal B}^{\bf P} (r) \} }{d r} \leq \frac{ d\; {\rm
vol} \{ {\cal B} (r) \} }{d r}$.  Hence inequality
$\frac{d \; {\rm vol}\{ {\cal B}^{\bf P} (r) \} }{d r} \leq
\Bbb{P}(r)\;\frac{ d\; {\rm vol} \{ {\cal B} (r) \} }{d r}$ is not hard to guarantee
 in the range of $r$ such that $\Bbb{P}(r)$ is close to $1$.
It follows from equation~(\ref{der}) that
$\frac{d \Bbb{P}(r)}{d r}
\leq 0$ can be readily satisfied in the case of small risk parameter $\epsilon$.

\bigskip

It is interesting to note
that our randomized algorithms for estimating $\rho(\epsilon)$
 presented in the sequel rely only on the following assumption:
\be
\Bbb{P} (r) < 1 - \epsilon  \;\;\forall r \in ( \rho(\epsilon) , 2^\kappa] \bigcup \{2^i: \kappa  < i \leq 0   \}
\label{simpler}
\ee
where $\kappa = \lceil \log_2 \rho(\epsilon) \rceil$.  It can be seen that condition~(\ref{simpler})
is even weaker than the separable assumption.

When condition~(\ref{simpler}) is guaranteed, an interval which includes
$\rho(\epsilon)$
can be readily found by starting from uncertainty radius $r=1$ and
then successively doubling $r$ or cutting $r$ in half based on the
comparison of $\Bbb{P}(r)$ with $1-\epsilon$.
Moreover, bisection method can be employed to refine the estimate for
$\rho(\epsilon)$.
Of course, the success of such methods depends
on the reliable and efficient comparison of $\Bbb{P}(r)$ with $1-\epsilon$ based on Monte
Carlo method.
In the following subsection, we illustrate a fast method of comparison.

\subsection{Probabilistic Comparison}

In general,  $\Bbb{P}_X$ can
only be estimated by a Monte Carlo method.
The conventional method is to directly compare $\frac{K}{N}$ with
$1 - \epsilon$ where $K$ is the number
of successful trials during $N$ i.i.d. sampling experiments.
  There are three problems with the conventional method.
  First,  the comparison of $\frac{K}{N}$ with $1-\epsilon$ can be very
misleading.
  Second, the sample size $N$ is required to be very large to obtain a
reliable comparison.
  Third, we don't know how reliable the comparison is.
  In this subsection, we present a new approach which
   allows for a reliable comparison with many fewer samples.
   The key idea is to compare binomial confidence limits with
   the fixed probability $1 - \epsilon$ and hence reliable judgement can be
made in advance.

\begin{description}

\item Function name: {\bf Probabilistic-Comparison}.

\item Input: Risk parameter $\epsilon$ and confidence parameter $\delta$.

\item Output: $d =$ {\bf Probabilistic-Comparison} $(\delta, \epsilon)$.

\item \underline{Step $1$.} Let $d \leftarrow 0$.

\item \underline{Step $2$.}  While $d=0$ do the following:

\begin{itemize}

\item Sample $X$.

\item Update $N$ and $K$.

\item Compute lower confidence limit $L$ and
and upper confidence limit $U$ by Theorem 1.

\item If $U < 1 - \epsilon$ then let $d \leftarrow -1$.
If $L > 1 - \epsilon$ then let $d \leftarrow 1$.

\end{itemize}

\end{description}

The confidence parameter $\delta$ is used to control the reliability
of the comparison.  A typical value is $\delta = 0.01$. The
implication of output is as follows: $d = 1$ indicates that
$\Bbb{P}_X
 > 1 - \epsilon$ is true with high confidence;
$d = - 1$ indicates that $\Bbb{P}_X < 1 - \epsilon$ is true with
high confidence.

\begin{figure}
\center

\includegraphics[height=8cm]{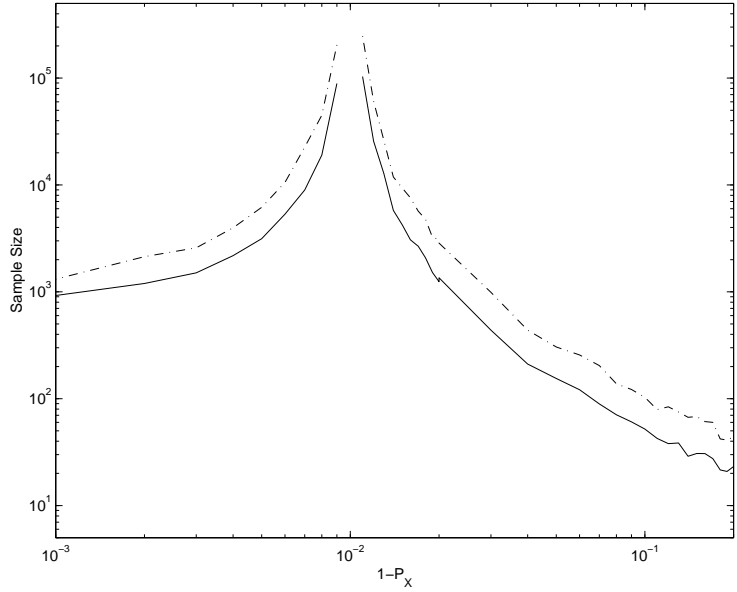}

\caption{Complexity of Probabilistic Comparison.  The horizontal
axis represents $1 - \Bbb{P}_X$.
 The vertical axis represents sample size.
The solid line and the dash-dot line respectively show the average
sample size and the $95\%$-quantile of the sample size. }
\label{fig48}       
\end{figure}

Obviously, the sample size is random in nature.  For $\epsilon =
\delta = 0.01$, we simulated the Probabilistic Comparison Algorithm
identically and independently $100$ times for different values of
$\Bbb{P}_X$. We observe that, for each value of $\Bbb{P}_X$, the
Probabilistic Comparison Algorithm makes correct judgement among all
$100$ simulations. Figure~\ref{fig48} shows the average sample size
and the $95\%$-quantile of the sample size estimated from our
simulation.  It can be seen from the figure that, as long as
$\Bbb{P}_X$ is not very close to $1 -\epsilon$, the Probabilistic
Comparison Algorithm can make a reliable comparison with a small
sample size.

\subsection{Computing Initial Interval}

Under the separable assumption, an interval $[a,b]$
which includes $\rho(\epsilon)$ can be quickly determined by the following
algorithm.

\begin{description}

\item Function name: {\bf Initial}.

\item Input: Risk parameter $\epsilon$ and confidence parameter $\delta$.

\item Output:
$[a,b] = {\bf Initial}(\delta,\epsilon)$.

\item \underline{Step $1$.} Let $r \leftarrow 1$.
Apply Probabilistic-Comparison algorithm to compare $\Bbb{P}(1)$ with $1 -
\epsilon$.  Let
the outcome be $d_1$.

\item \underline{Step $2$.}  If $d_1=1$ then let $d  \leftarrow  d_1$ and do
the following:

\begin{itemize}

\item While $d=1$ do the following:
       \begin{itemize}
       \item Let $r  \leftarrow  2r$.  Apply Probabilistic-Comparison
algorithm to compare
$\Bbb{P}(r)$ with $1 - \epsilon$.  Let
the outcome be $d$.
       \end{itemize}

\item Let $a  \leftarrow \frac{r}{2}$ and $b \leftarrow r$.
\end{itemize}

\item \underline{Step $3$.}  If $d_1=-1$ then let $d  \leftarrow  d_1$ and
do the following:

\begin{itemize}

\item While $d=-1$ do the following:
       \begin{itemize}
       \item Let $r \leftarrow  \frac{r}{2}$.  Apply
Probabilistic-Comparison algorithm to compare
$\Bbb{P}(r)$ with $1 - \epsilon$.  Let the outcome be $d$.

       \end{itemize}

\item Let $a  \leftarrow r$ and $b \leftarrow 2r$.
\end{itemize}

\end{description}


\subsection{Probabilistic Bisection}

Once an initial interval $[a,b]$ is obtained, an estimate $\widehat{R}$ for
the probabilistic robustness margin $\rho(\epsilon)$ can be efficiently
computed as follows.

\begin{description}

\item Function name: {\bf Bisection}.

\item Input: Risk parameter $\epsilon$, confidence parameter $\delta$,
initial interval $[a,b]$, relative tolerance $\gamma$.

\item Output: $[\widehat{R}] ={\bf Bisection}(\gamma, a,b,\delta,\epsilon)$.

\item \underline{Step $1$.} Input
$\gamma,\;a,\;b,\;\delta,\;\epsilon$.

\item \underline{Step $2$.} While $b-a  > \gamma  a$
do the following:

      \begin{itemize}

      \item
Let $r \leftarrow \frac{a+b}{2}$.   Apply Probabilistic-Comparison algorithm
to compare
$\Bbb{P}(r)$ with $1 - \epsilon$.  Let
the outcome be $d$.

      \item If $d=-1$ then let $b \leftarrow r$, else let $a \leftarrow r$.

      \end{itemize}

\item \underline{Step $3$.} Return $ \widehat{R} = b$ (Note: this is
actually a soft upper bound).

\end{description}


It should be noted that when applying bisection algorithm to refine the initial interval $[a, \; b]$,
the execution of the algorithm may take very long time if $\Bbb{P} (r) \approx \epsilon$.
However, such chance is almost $0$.  This problem can be fixed by the following methods.

\begin{itemize}

\item We can limit the maximum number of simulations to a number $M$.
When conducting simulation at radius $r$,  simulation results
can be saved for uncertainty radius $\frac{a + r} {2}$
where $a$ is the lower bound of the current interval with middle point $r$
(See the next section for the idea of sample reuse).
After the number of simulations for uncertainty radius $r$ exceeds $M$,
the simulation is switched for uncertainty radius $\frac{a + r} {2}$.

\item In application, one might want to construct the robustness degradation curve in
the backward direction by starting from an upper bound of $\rho(\epsilon)$.
In that situation, we don't need to compute a very tight
interval for $ \rho (\epsilon)$ and hence the chance of $\Bbb{P} (r) \approx \epsilon$ is even smaller.

\end{itemize}

\section{Constructing Robustness Degradation Curve}

We shall develop efficient randomized algorithms for constructing robustness
degradation
curve, which provide more insight for the robustness of the uncertain system
than probabilistic robustness margin.
First we recall the Sample Reuse algorithm \cite{C3} for
constructing robustness degradation curve for a given range of uncertainty
radius.

\begin{description}

\item \underline{{\bf Sample Reuse Algorithm}}

\item Input: Sample size $N$, confidence parameter $\delta \in (0,1)$,
uncertainty radius interval $[a, b]$, number of uncertainty radii $l$.

\item Output: Proportion estimate $\widehat{\Bbb{P}}_i$ and the related
confidence interval
for $r_i = b - \frac{(b-a)(i-1)}{l-1}, \;\; i=1, 2, \cdots, l$.  In the
following, $m_{i1}$ denotes the number of sampling experiments
conducted at $r_i$ and
$m_{i2}$ denotes the number of observations guaranteeing ${\bf P}$
during the $m_{i1}$ sampling experiments.

\item \underline{Step $1$ (Initialization).}
Let $M =[m_{ij}]_{l \times 2}$ be a zero matrix.

\item \underline{Step $2$ (Backward Iteration).} For $i=1$ to $i=l$ do the
following:

      \begin{itemize}

      \item Let $r  \leftarrow  r_i$.

      \item While $m_{i1} < N$ do the following:

              \begin{itemize}

                \item Generate uniform
                  sample $q$ from ${\cal B}(r)$. Evaluate the robustness requirement
                  {\bf P} for $q$.

                \item Let $m_{s1} \leftarrow m_{s1}+1$ for any $s$ such that
$ r_s \geq \ell(q)$.

                \item If robustness requirement {\bf P} is satisfied
                      for $q$ then let $m_{s2} \leftarrow m_{s2}+1$ for
any $s$ such that $ r_s \geq \ell(q)$.

              \end{itemize}

      \item  Let $\widehat{\Bbb{P}}_i  \leftarrow \frac{m_{i2}}{N}$
      and construct the confidence interval of confidence level $100(1-\delta) \%$ by Theorem $1$.

      \end{itemize}

\end{description}

\bigskip
It should be noted that the idea of the Sample Reuse Algorithm is not simply
a save of sample generation.
It is actually a backward iterative mechanism.
In the algorithm, the most important save of computation is usually
the evaluation of the complex robustness requirements ${\bf P}$ (See
\cite{C3} for details).

Now we introduce the global strategy for constructing robustness degradation
curve.  The idea is to successively apply the Sample Reuse Algorithm
for a sequence of intervals of uncertainty radius. Each time the size of
interval is reduced by half.
The lower bound of the current interval is defined to be the upper bound of
the next consecutive interval.  The algorithm is
terminated once the robustness requirement ${\bf P}$ is guaranteed for all $N$ samples of
an uncertainty set of which the radius is taken as the lower bound of an interval of uncertainty radius.
More precisely, the procedure is presented as follows.

\begin{description}

\item \underline{{\bf Global Strategy}}

\item Input: Sample size $N$, risk parameter $\epsilon$ and confidence
parameter $\delta \in (0,1)$.

\item Output: Proportion estimate $\widehat{\Bbb{P}}_i$ and the related
confidence interval.

\item \underline{Step $1$.} Compute an estimate $\widehat{R}$ for
probabilistic robustness margin $\rho(\epsilon)$.

\item \underline{Step $2$.} Let $STOP \leftarrow 0$. Let $a \leftarrow
\frac{\widehat{R}} {2}$ and $b \leftarrow \widehat{R}$.

\item \underline{Step $3$ (Backward Iteration).} While $STOP =0$ do the
following:

      \begin{itemize}

      \item Apply Sample Reuse Algorithm to construct robustness degradation
curve for uncertainty radius interval $[a, b]$.

      \item  If the robustness property ${\bf P}$ is guaranteed for all $N$
samples of uncertainty set ${\cal B}(r)$ with radius $r = a$
      then let $STOP \leftarrow 1$, otherwise let $b \leftarrow a$ and $a
\leftarrow \frac{b}{2}$.

      \end{itemize}

\end{description}

\bigskip

For given risk parameter $\epsilon$ and confidence parameter $\delta$,
the sample size is chosen as
\be
N > \frac{ 2 (1 - \epsilon + \frac{\alpha \epsilon }{3} )
(1 - \frac{\alpha }{3} ) \ln \frac{2}{\delta} } { \alpha^2 \epsilon }
\label{size}
\ee
with $\alpha \in (0,1)$.  It follows from Massart's inequality
\cite{Massart} that such a sample size ensures
\[
\Pr \left\{ \left| \Bbb{P}_X - \frac{K}{N}  \right| < \alpha \epsilon
\right\} > 1 - \delta
\]
with $\Bbb{P}_X = 1 - \epsilon$ (See also Lemma~\ref{Massart} in Appendix).
It should be noted that Massart's inequality
is less conservative than the Chernoff bounds in both multiplicative and
additive forms.

\bigskip

We would like to remark that the algorithms we propose for estimating the
probabilistic robustness margin and
constructing robustness degradation curve are susceptible of further improvement.
For example, the idea of sample reuse is not employed in computing the
initial interval
and in the probabilistic bisection algorithm.  Moreover, in constructing the
robustness degradation curve,
the Sample Reuse Algorithm is independently applied for each interval of
uncertainty radius.
Actually, the simulation results can be saved for the successive intervals.

We would also like to note that, for a practitioner,
computing the probabilistic robustness margin might be sufficient for understanding
 the system robustness.  However,  when more insight about the system robustness is expected,
 the techniques introduced in Section $5$ can be employed.
 Of course, the price is more computational effort.

\section{Illustrative Examples}

In this section we demonstrate through examples the power of our
randomized algorithms
in solving a wide variety of complicated robustness analysis problems which
are
not tractable in the classical deterministic framework.

\begin{figure}[htbp]
\centerline{\psfig{figure=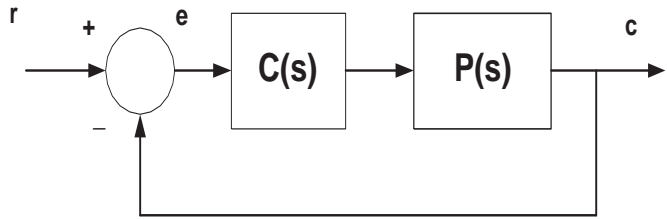, height=1.3in, width=4in
}} \caption{Uncertain System }
\label{fig_08}
\end{figure}

We consider a system which has been studied in \cite{GS} by a deterministic
approach.
The system is as shown in Figure~\ref{fig_08}.

The compensator is $C(s)=\frac{s+2}{s+10}$ and the plant is
$P(s)=\frac{800(1+0.1\delta_1)}{s(s+4+0.2\delta_2)(s+6+0.3\delta_3)}$
with parametric uncertainty $\Delta=[\delta_1,\delta_2,\delta_3]^{\rm T}$.
The nominal system is stable.
The closed-loop roots of the nominal system are:
\[
z_1= -15.9178,  \;\;  z_2 = -1.8309,  \;\; z_3 = -1.1256 + 7.3234i,  \;\;
z_4 = -1.1256 - 7.3234i.
\]
The peak value, rise time, settling time of step response of the
nominal system,  are respectively,
$P_{peak}^0 =1.47, \;\;  t_r^0 = 0.185, \;\; t_s^0 = 3.175$.

We first consider the robust ${\cal D}$-stability of the system.
The robustness
requirement ${\bf P}$ is defined as
${\cal D}$-stability with the domain of poles defined as: Real part $ <
-1.5$, or fall within one of the
two disks centered at $z_3$
and $z_4$ with radius $0.3$.  The uncertainty set is defined as the polytope
\[
{\cal B}_H (r):=\left\{ r \Delta +(1-r)\;\frac{ \sum_{i=1}^4 \Delta^i} {4} :
\Delta \in {\rm conv}\{\Delta^1,\Delta^2,\Delta^3,\Delta^4\}  \right\}
\]
where `conv' denotes the convex hull of $\Delta^i=[ \frac{1}{2}
\sin(\frac{2i-1}{3}\pi),\;
\frac{1}{2} \cos(\frac{2i-1}{3}\pi),\;
-\frac{\sqrt{3}}{2}]^{\rm T}$ for $i=1,2,3$ and
$\Delta^4=[
0,\;0,\;1]^{\rm T}$.

\begin{figure}[htbp]
\centerline{\psfig{figure=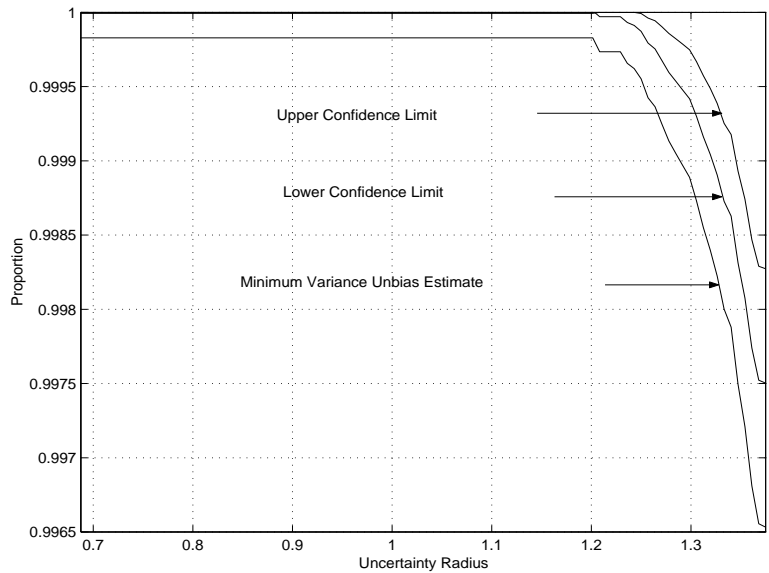, height=3in, width=4in
}} \caption{ Robustness Degradation Curve}
\label{fig_000}
\end{figure}

Obviously, there exists no effective method for computing the
 deterministic robustness margin in the literature.  However,
 our randomized algorithms can efficiently construct the robustness degradation curve.  See
 Figure~\ref{fig_000}.

 In this example, the risk parameter is {\it a priori} specified as $\epsilon
= 0.001$.  The procedure for estimating the probabilistic robustness margin is
explained as follows.  Let $N$ denote the sample size which is random in nature.  Let $K$
denote the number of successful trials among $N$ i.i.d. sampling experiments as defined
in Section $2$ and Subsection $3.1$ (i.e., a successful trial is equivalent to
an observation that the robustness requirement is guaranteed).
Let confidence parameter $\delta = 0.01$ and choose tolerance $\gamma = 0.05$.
Staring from $r=1$, after $N = 7060$ simulations we obtain $K= 7060$, the
probabilistic comparison
algorithm determined that
$\Bbb{P}(1) > 1 - \epsilon$ since the lower confidence limit $L > 1 -
\epsilon$.
The simulation is thus switched to uncertainty radius $r=2$. After $N = 65$
times of simulation,
it is found that $K = 61$.  The probabilistic
comparison algorithm detected that $\Bbb{P}(2) < 1- \epsilon$ because the
upper confidence limit $U < 1- \epsilon$.
So, initial interval $[1,2]$ is readily obtained.  Now the probabilistic
bisection algorithm is invoked.
Staring with the middle point of the initial interval (i.e., $r= \frac{1+2}{2} = \frac{3}{2}$),
after $N = 613$ times of simulations, it is found that $K = 607$,
the probabilistic comparison algorithm concluded that $\Bbb{P}(\frac{3}{2}) < 1 -
\epsilon$
since the upper confidence limit $U < 1 - \epsilon$.
Thus simulation is moved to $r = \frac{1 + \frac{3}{2}}{2} = \frac{5}{4}$.
It is found that $K = 9330$ among $N = 9331$ times of simulations.
Hence, the probabilistic comparison algorithm determined that $\Bbb{P}(\frac{5}{4})
 > 1 - \epsilon$
since the lower confidence limit $L > 1 - \epsilon$.
Now the simulation is performed at $r = \frac{\frac{5}{4} + \frac{3}{2}}{2} = \frac{11}{8}$.
After $N = 6653$ simulations,
it is discovered that $K = 6636$. The probabilistic comparison algorithm
judged that $\Bbb{P}(\frac{11}{8}) < 1 - \epsilon$
based on calculation that the upper confidence limit $U < 1- \epsilon$.
At this point the interval is $[\frac{5}{4}, \frac{11}{8}]$
and the bisection is terminated since tolerance condition $b-a \leq \gamma
a$ is satisfied.
The evolution of intervals produced by the probabilistic bisection algorithm
is as follows:
\[
[1,2] \;\; \longrightarrow \;\; \left[1, \frac{3}{2}\right] \;\; \longrightarrow
\;\;
\left[ \frac{5}{4}, \frac{3}{2} \right] \;\; \longrightarrow \;\; \left[
\frac{5}{4}, \frac{11}{8} \right].
\]

Now we have obtained an interval $[\frac{5}{4}, \frac{11}{8}]$ which includes
$\rho(0.001)$, so the Sample Reuse
Algorithm can be employed to construct robustness degradation curve.
In this example, the number of uncertainty radii is $l = 100$ and the
confidence parameter is chosen as $\delta = 0.001$.
A constant sample size is computed by formula~(\ref{size}) with $\alpha = 0.5$ as
\[
N = 50,631.
\]
The interval from which we start constructing robustness degradation curve is
$[\frac{11}{16},\frac{11}{8}]$.
It is determined that $K = N = 50,632$ at uncertainty radius $r = \frac{11}{16}$.
Therefore,
the Sample Reuse Algorithm is invoked only once and the overall algorithm is
terminated
(If $K \neq N$ for $r = \frac{11}{16}$, then the next interval would be $[\frac{11}{32}, \frac{11}{16}]$).
Although $\Bbb{P}(r)$ is evaluated for $l= 100$ uncertainty radii with the
same sample size $N$,
the total number of simulations is only $153,358 << N l = 100 N$.
To provide an evaluation of accuracy for all estimates of $\Bbb{P}(r)$, confidence limits
are computed by Theorem $1$.

We now apply our algorithms to
 a robustness problem with time specifications. Specifically, the robustness requirement ${\bf
P}$ is : Stability, and rise time
$t_r < 135 \% \; t_r^0  = 0.25$,
settling time $t_s < 110 \% \; t_s^0 = 3.5$,
overshoot $P_{peak} < 116 \% \; P_{peak}^0= 1.7$.  The uncertainty set is
${\cal B}_{\infty} (r) :=
\{ \Delta:\;||\Delta||_{\infty} \leq r \}$.

\begin{figure}[htbp]
\centerline{\psfig{figure=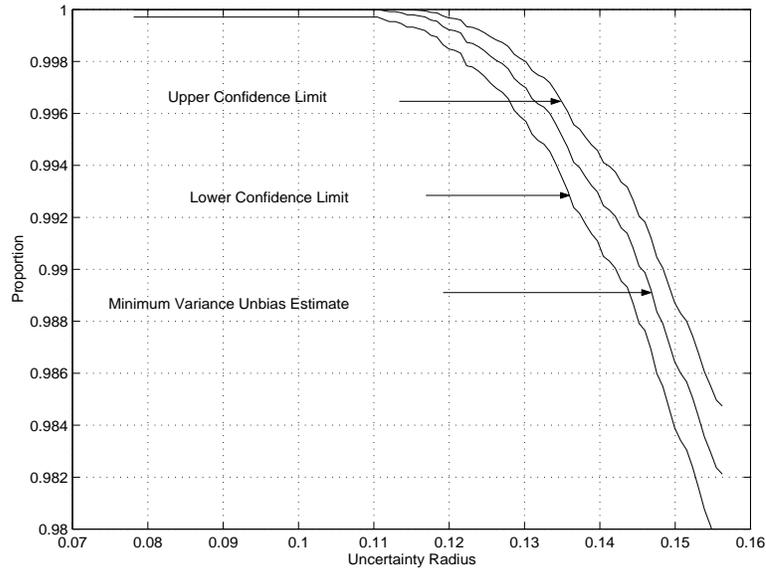, height=3in, width=4in
}} \caption{ Robustness Degradation Curve}
\label{fig_0000}
\end{figure}

In this case, the risk parameter is {\it a priori} specified as $\epsilon
= 0.01$.  It is well known that, for this type of problem, there exists no effective method for computing the
 deterministic robustness margin in the literature.  However,
 our randomized algorithms can efficiently construct the robustness degradation curve.  See Figure~\ref{fig_0000}.

 We choose $\gamma = 0.25$ and $\delta = 0.01$ for estimating $\rho(0.01)$.
 Starting from uncertainty radius $r=1$, the initial interval is easily found as
 $[\frac{1}{8}, \frac{3}{16}]$ through the following interval evolution:
\[
\left[ \frac{1}{2},1 \right] \;\; \longrightarrow \;\; \left[\frac{1}{4}, \frac{1}{2}\right] \;\; \longrightarrow
\;\;
\left[ \frac{1}{8}, \frac{1}{4} \right].
\]
The sequence of intervals generated by the probabilistic bisection algorithm
is as follows:
\[
\left[\frac{1}{8}, \frac{1}{4}\right] \;\; \longrightarrow
\;\;
\left[ \frac{1}{8}, \frac{3}{16} \right] \;\; \longrightarrow
\;\;
\left[ \frac{1}{8}, \frac{5}{32} \right].
\]
So, we obtained an estimate for the probabilistic robustness margin $\rho(0.01)$ as $\frac{5}{32}$.
To construct robustness degradation curve, the number of uncertainty radii is
chosen as $l = 100$ and the
confidence parameter is chosen as $\delta = 0.01$.
A constant sample size is computed by formula~(\ref{size}) with $\alpha = 0.2$ as
\[
N = 24,495.
\]
The interval from which we start constructing robustness degradation curve is
$[\frac{5}{64},\frac{5}{32}]$.  We found that this is also the last interval of uncertainty radius because
it is determined that $K = N$ at uncertainty radius $r = \frac{5}{64}$.

\section{Conclusions}

In this paper, we have established efficient techniques which applies to
robustness analysis problems with arbitrary robustness requirements and
uncertainty bounding set.
The key mechanisms are probabilistic comparison, probabilistic bisection and
backward iteration.
Motivated by the crucial role of binomial confidence interval in reducing
the computational complexity, we have derived an explicit formula for
computing binomial confidence limits.
This formula overcomes the computational issue and inaccuracy of standard
methods.

\bigskip

\appendix

\sect{Proof of Theorem 1}

\bigskip

To show Theorem~\ref{Tape_Massart}, we need some preliminary results.
The following Lemma \ref{Massart} is due to Massart \cite{Massart}.

\begin{lemma} \label{Massart}
${\Pr} \left\{ \frac{K}{N} \geq \Bbb{P}_X + \epsilon \right\} \leq
\exp \left(  - \frac{N \epsilon^2} {2 ( \Bbb{P}_X +
\frac{\epsilon}{3} )\;( 1 - \Bbb{P}_X - \frac{\epsilon}{3} )}
\right)$ for all  $\epsilon \in (0, 1 - \Bbb{P}_X)$.
\end{lemma}

Of course, the above upper bound holds trivially for $\epsilon \geq 1 -
\Bbb{P}_X$.
Thus, Lemma~\ref{Massart} is actually true for any $\epsilon >0$.

\begin{lemma} \label{Massart2} ${\Pr} \left\{ \frac{K}{N} \leq \Bbb{P}_X -
\epsilon \right\} \leq
\exp \left(  - \frac{N \epsilon^2} {2 ( \Bbb{P}_X -
\frac{\epsilon}{3} )\;( 1 - \Bbb{P}_X + \frac{\epsilon}{3} )}
\right)$ for all  $\epsilon > 0$.
\end{lemma}

\begin{pf}

Define $Y = 1 - X$.  Then $\Bbb{P}_Y = 1 - \Bbb{P}_X$.  At the same time
when we are
conducting $N$ i.i.d. experiments for $X$, we are also conducting $N$ i.i.d.
experiments for $Y$.
Let the number of successful trials of the experiments for $Y$ be denoted as
$K_Y$.  Obviously, $K_Y = N - K$.
Applying Lemma~\ref{Massart} to $Y$, we have
$
{\Pr} \left\{ \frac{K_Y}{N} \geq \Bbb{P}_Y + \epsilon \right\} \leq
\exp \left(  - \frac{N \epsilon^2} {2 ( \Bbb{P}_Y +
\frac{\epsilon}{3} )\;( 1 - \Bbb{P}_Y - \frac{\epsilon}{3} )}
\right)$.
It follows that
$
{\Pr} \left\{ \frac{N - K}{N} \geq 1 - \Bbb{P}_X + \epsilon \right\} \leq
\exp \left(  - \frac{N \epsilon^2} {2 ( 1 - \Bbb{P}_X +
\frac{\epsilon}{3} )\;[ 1 - (1- \Bbb{P}_X) - \frac{\epsilon}{3} ]}
\right)$.
The proof is thus completed by observing that ${\Pr} \left\{ \frac{N - K}{N}
\geq 1 - \Bbb{P}_X + \epsilon \right\} =
{\Pr} \left\{ \frac{K}{N} \leq \Bbb{P}_X - \epsilon \right\}$.
\end{pf}

The following lemma can be found in \cite{CR}.
\begin{lemma} \label{decrease} $\sum_{j=0}^{k} {N \choose j}
   x^j (1- x)^{N-j}$ decreases monotonically with respect to $x \in (0,1)$
for
$k=0,1,\cdots,N$.
\end{lemma}

\begin{lemma} \label{lem_massart} $ \sum_{j=0}^{k} {N \choose j}
   x^j (1- x)^{N-j} \leq \exp \left( - \frac{ N ( x - \frac{k}{N} )^2 }
{ 2 \; ( \frac{2}{3} x + \frac{k}{3N} ) \; ( 1- \frac{2}{3} x - \frac{k}{3N}
) } \right)
\;\;\forall x \in (\frac{k}{N}, 1)$ for
$k=0,1,\cdots,N$.
\end{lemma}
{\bf Proof.}  Consider binomial random variable $X$ with parameter
$\Bbb{P}_X
>\frac{k}{N}$.  Let $K$ be the number
of successful trials
during $N$ i.i.d. sampling experiments.  Then
$\sum_{j=0}^{k} {N \choose j}
   \Bbb{P}_X^j (1- \Bbb{P}_X)^{N-j} = {\rm Pr} \{ K \leq k\}$.
   Note that ${\rm Pr} \{ K \leq k\} = {\rm Pr} \left\{ \frac{K}{N} \leq
\Bbb{P}_X - \left( \Bbb{P}_X -
\frac{k}{N} \right) \right\}$. Applying Lemma \ref{Massart2} with
$\epsilon = \Bbb{P}_X - \frac{k}{N} > 0$, we have
\begin{eqnarray*}
\sum_{j=0}^{k} {N \choose j}
   \Bbb{P}_X^j (1- \Bbb{P}_X)^{N-j}
& \leq & \exp \left(  - \frac{N ( \Bbb{P}_X - \frac{k}{N} )^2} {2 (
\Bbb{P}_X -
\frac{\Bbb{P}_X - \frac{k}{N}}{3} )\;( 1 - \Bbb{P}_X + \frac{\Bbb{P}_X -
\frac{k}{N}}{3} )}
\right)\\
& = & \exp \left( - \frac{ N ( \Bbb{P}_X  - \frac{k}{N} )^2 }
{ 2 \; ( \frac{2}{3} \Bbb{P}_X  + \frac{k}{3N} ) \;
( 1- \frac{2}{3} \Bbb{P}_X  - \frac{k}{3N} ) } \right).
\end{eqnarray*}
Since the argument holds for arbitrary binomial random variable $X$ with
$\Bbb{P}_X
>\frac{k}{N}$, thus the proof of the lemma is completed.
$\;\;\;\;\square$

\begin{lemma}\label{lem_Massart_B}
$\sum_{j=0}^{k-1} {N \choose j}
   x^j (1- x)^{N-j}
   \geq 1 - \exp \left( - \frac{ N ( x - \frac{k}{N} )^2 }
{ 2 \; ( \frac{2}{3} x + \frac{k}{3N} ) \; ( 1- \frac{2}{3} x - \frac{k}{3N}
) } \right)
\;\;
\forall x \in (0,\frac{k}{N})$ for
$k=1,\cdots,N$.
\end{lemma}
{\bf Proof.}  Consider binomial random variable $X$ with parameter
$ \Bbb{P}_X < \frac{k}{N}$.  Let $K$ be the number
of successful trials
during $N$ i.i.d. sampling experiments.
Then
\[
\sum_{j=0}^{k-1} {N \choose j}
   \Bbb{P}_X^j (1- \Bbb{P}_X)^{N-j}
=  {\rm Pr} \{K < k\} =
{\rm Pr} \left\{ \frac{K}{N} < \Bbb{P}_X  +
(\frac{k}{N} - \Bbb{P}_X) \right\}.
\]
Applying Lemma \ref{Massart} with
$\epsilon = \frac{k}{N} - \Bbb{P}_X > 0$, we have that
\begin{eqnarray*}
\sum_{j=0}^{k-1} {N \choose j}
   \Bbb{P}_X^j (1- \Bbb{P}_X)^{N-j}
& \geq & 1 - \exp \left(  - \frac{N (\frac{k}{N} - \Bbb{P}_X )^2}
{2 ( \Bbb{P}_X +
\frac{ \frac{k}{N} - \Bbb{P}_X }{3} )\;
( 1 - \Bbb{P}_X - \frac{ \frac{k}{N} - \Bbb{P}_X }{3} )}
\right)\\
& = & 1 - \exp \left( - \frac{ N ( \Bbb{P}_X  - \frac{k}{N} )^2 }
{ 2 \; ( \frac{2}{3} \Bbb{P}_X  + \frac{k}{3N} ) \;
( 1- \frac{2}{3} \Bbb{P}_X  - \frac{k}{3N} ) } \right).
\end{eqnarray*}
Since the argument holds for arbitrary binomial random variable $X$
with $\Bbb{P}_X < \frac{k}{N}$,
thus the proof of the lemma is completed.
$\;\;\;\;\square$

\begin{lemma}\label{lem8} Let $0 \leq k \leq N$. Then
$L_{N,k,\delta} < U_{N,k,\delta}$.
\end{lemma}

{\bf Proof.}  Obviously, the lemma is true for $k=0,N$.
We consider the case
that $1 \leq k \leq N-1$.
Define ${\cal S} (N,k,x) := \sum_{j=0}^{k} {N \choose j}
x^j (1- x)^{N-j}$ for $x \in (0,1)$.  Notice that
${\cal S} (N,k,\overline{p})=
{\cal S} (N,k-1,\overline{p})+{N \choose
k}
   \overline{p}^k (1-\overline{p})^{N-k}=
   \frac{\delta}{2}$.  Thus
\[
{\cal S} (N,k-1,\underline{p})-{\cal S} (N,k-1,\overline{p})=
   1-\frac{\delta}{2}-\left[\frac{\delta}{2}-{N \choose k}
   \overline{p}^k (1-\overline{p})^{N-k}\right].
   \]
Notice that $\delta \in (0,1)$ and that $\overline{p} \in (0,1)$,
we have that
\[
{\cal S} (N,k-1,\underline{p})-{\cal S} (N,k-1,\overline{p})
= 1-\delta+
   {N \choose k} \overline{p}^k (1-\overline{p})^{N-k} > 0.
\]
By Lemma~\ref{decrease}, ${\cal S} (N,k-1,x)$
decreases monotonically with respect to $x$,
we have that $\underline{p} <
\overline{p}$.   $\;\;\;\;\square$

\bigskip

We are now in the position to prove Theorem~\ref{Tape_Massart}.  For national simplicity, let
\[
p = U_{N,k,\delta}, \;\;\; q = {\cal U} (N,k,\delta).
\]
It can be easily verified that $p \leq q$ for $k =
0, \; N$.
We need to show that $p \leq q$ for $0 < k < N$.
Straightforward computation shows that $q$ is the only root of
equation
\[
\exp \left( - \frac{ N ( x - \frac{k}{N} )^2 }
{ 2 \; ( \frac{2}{3} x + \frac{k}{3N} ) \; ( 1- \frac{2}{3} x - \frac{k}{3N}
) } \right)
= \frac{\delta}{2}
\]
with
respect to $x \in (\frac{k}{N},\infty)$.  There are two cases: $q
\geq 1$ and $q < 1$.
If $q \geq 1$ then $p \leq q$
is trivially true.  We only need to consider the case that $\frac{k}{N} <
q < 1$.
In this case, it follows from Lemma~\ref{lem_massart}
that
\[
\sum_{j=0}^{k} {N \choose j}
   {q}^j (1- q)^{N-j} \leq \exp \left( - \frac{ N (
q - \frac{k}{N} )^2 }
{ 2 \; ( \frac{2}{3} q + \frac{k}{3N} ) \; ( 1- \frac{2}{3} q  - \frac{k}{3N} ) } \right)
= \frac{\delta}{2}.
\]
Recall that
\[
\sum_{j=0}^{k} {N \choose j}
   p^j (1- p)^{N-j} = \frac{\delta}{2},
   \]
we have
\[
\sum_{j=0}^{k} {N \choose j}
   p^j (1-p)^{N-j} \geq \sum_{j=0}^{k} {N \choose
j}
   {q}^j (1- q)^{N-j}.
   \]
   Therefore, by Lemma~\ref{decrease}, we have that $p \leq q$ for $0 < k < N$.
Thus, we have shown that $U_{N,k,\delta} \leq {\cal U}_{N,k,\delta}$ for all $k$.

Similarly, by Lemma~\ref{lem_Massart_B} and Lemma~\ref{decrease} we can show
that
$L_{N,k,\delta} \geq  {\cal L}(N,k,\delta)$.  By Lemma~\ref{lem8}, we have
${\cal L}(N,k,\delta) < L_{N, k, \delta} < U_{N, k, \delta} < {\cal U}(N,k,\delta)$.
Finally, the proof of Theorem~\ref{Tape_Massart} is
completed by invoking the probabilistic
implication of the Clopper-Pearson confidence interval.

\end{document}